\documentclass[11pt]{amsart}
\usepackage[parfill]{parskip} 
\usepackage{amsmath, amsxtra, amsthm, amssymb,eucal,dsfont}
\usepackage[all]{xy}
\usepackage{color}
\usepackage[colorlinks=true, pdfstartview=FitV, linkcolor=blue, citecolor=blue, urlcolor=blue]{hyperref}

\usepackage{graphicx}

\numberwithin{equation}{section}
\swapnumbers

\newcommand{\nc}{\newcommand}
\nc{\rc}{\renewcommand}

\rc{\b}{\mathbb}
\rc{\c}{\mathcal}
\nc{\on}{\operatorname}
\nc{\tn}{\textnormal}

\nc{\bA}{\b A}
\nc{\bB}{\b B}
\nc{\bC}{\b C}
\nc{\bD}{\b D}
\nc{\bE}{\b E}
\nc{\bF}{\b F}
\nc{\bG}{\b G}
\nc{\bH}{\b H}
\nc{\bI}{\b I}
\nc{\bJ}{\b J}
\nc{\bK}{\b K}
\nc{\bL}{\b L}
\nc{\bM}{\b M}
\nc{\bN}{\b N}
\nc{\bO}{\b O}
\nc{\bP}{\b P}
\nc{\bQ}{\b Q}
\nc{\bR}{\b R}
\nc{\bS}{\b S}
\nc{\bT}{\b T}
\nc{\bU}{\b U}
\nc{\bV}{\b V}
\nc{\bW}{\b W}
\nc{\bX}{\b X}
\nc{\bY}{\b Y}
\nc{\bZ}{\b Z}

\nc{\cA}{\c A}
\nc{\cB}{\c B}
\nc{\cC}{\c C}
\nc{\cD}{\c D}
\nc{\cE}{\c E}
\nc{\cF}{\c F}
\nc{\cG}{\c G}
\nc{\cH}{\c H}
\nc{\cI}{\c I}
\nc{\cJ}{\c J}
\nc{\cK}{\c K}
\nc{\cL}{\c L}
\nc{\cM}{\c M}
\nc{\cN}{\c N}
\nc{\cO}{\c O}
\nc{\cP}{\c P}
\nc{\cQ}{\c Q}
\nc{\cR}{\c R}
\nc{\cS}{\c S}
\nc{\cT}{\c T}
\nc{\cU}{\c U}
\nc{\cV}{\c V}
\nc{\cW}{\c W}
\nc{\cX}{\c X}
\nc{\cY}{\c Y}
\nc{\cZ}{\c Z}

\nc{\fA}{{\mathfrak A}}
\nc{\fB}{{\mathfrak B}}
\nc{\fC}{{\mathfrak C}}
\nc{\fD}{{\mathfrak D}}
\nc{\fE}{{\mathfrak E}}
\nc{\fF}{{\mathfrak F}}
\nc{\fG}{{\mathfrak G}}
\nc{\fH}{{\mathfrak H}}
\nc{\fI}{{\mathfrak I}}
\nc{\fJ}{{\mathfrak J}}
\nc{\fK}{{\mathfrak K}}
\nc{\fL}{{\mathfrak L}}
\nc{\fM}{{\mathfrak M}}
\nc{\fN}{{\mathfrak N}}
\nc{\fO}{{\mathfrak O}}
\nc{\fP}{{\mathfrak P}}
\nc{\fQ}{{\mathfrak Q}}
\nc{\fR}{{\mathfrak R}}
\nc{\fS}{{\mathfrak S}}
\nc{\fT}{{\mathfrak T}}
\nc{\fU}{{\mathfrak U}}
\nc{\fV}{{\mathfrak V}}
\nc{\fW}{{\mathfrak W}}
\nc{\fZ}{{\mathfrak Z}}
\nc{\fX}{{\mathfrak X}}
\nc{\fY}{{\mathfrak Y}}
\nc{\fa}{{\mathfrak a}}
\nc{\fb}{{\mathfrak b}}
\nc{\fc}{{\mathfrak c}}
\nc{\fd}{{\mathfrak d}}
\nc{\fe}{{\mathfrak e}}
\nc{\ff}{{\mathfrak f}}
\nc{\fg}{{\mathfrak g}}
\nc{\fh}{{\mathfrak h}}
\nc{\fiI}{{\mathfrak i}}  
\nc{\ffi}{{\mathfrak i}}  
\nc{\fj}{{\mathfrak j}}
\nc{\fk}{{\mathfrak k}}
\nc{\fl}{{\mathfrak{l}}}
\nc{\fm}{{\mathfrak m}}
\nc{\fn}{{\mathfrak n}}
\nc{\fo}{{\mathfrak o}}
\nc{\fp}{{\mathfrak p}}
\nc{\fq}{{\mathfrak q}}
\nc{\fr}{{\mathfrak r}}
\nc{\fs}{{\mathfrak s}}
\nc{\ft}{{\mathfrak t}}
\nc{\fu}{{\mathfrak u}}
\nc{\fv}{{\mathfrak v}}
\nc{\fw}{{\mathfrak w}}
\nc{\fz}{{\mathfrak z}}
\nc{\fx}{{\mathfrak x}}
\nc{\fy}{{\mathfrak y}}

\nc{\al}{{\alpha }}
\nc{\be}{{\beta }}
\nc{\ga}{{\gamma }}
\nc{\de}{{\delta }}
\nc{\del}{{\partial }}
\nc{\vep}{{\varepsilon }}
\nc{\ep}{{\epsilon }}

\nc{\ze}{{\zeta }}
\nc{\et}{{\eta }}
\rc{\th}{{\theta }}

\nc{\vth}{{\vartheta }}

\nc{\io}{{\iota }}
\nc{\ka}{{\kappa }}
\nc{\la}{{\lambda }}
\nc{\vrho}{{\varrho}}
\nc{\si}{{\sigma }}
\nc{\ups}{{\upsilon }}
\nc{\vphi}{{\varphi }}
\nc{\om}{{\omega }}

\nc{\Ga}{{\Gamma }}
\nc{\De}{{\Delta }}
\nc{\nab}{{\nabla}}
\nc{\Th}{{\Theta }}
\nc{\La}{{\Lambda }}
\nc{\Si}{{\Sigma }}
\nc{\Ups}{{\Upsilon }}
\nc{\Om}{{\Omega }}

\nc{\Spec}{\on{Spec}}
\nc{\id}{\on{id}}
\nc{\inv}{ ^{-1}}
\nc{\su}{\subset}
\nc{\ot}{\otimes}
\nc{\un}{\underline}
\nc{\ov}{\overline}
\nc{\uu}{\mathds{1}}

\nc{\shom}{\cH om}
\nc{\Gr}{\cG\tn{r}}

\nc{\se}{\section}
\nc{\sse}{\subsection}
\nc{\ssse}{\subsubsection}

\nc{\md}{\!\!\!\mod}
\nc{\wt}{\widetilde}
\nc{\slc}{\fs \fl_2\bC}
\nc{\gr}{\tn{gr}}
\nc{\lan}{\langle}
\nc{\ran}{\rangle}
\nc{\ord}{\tn{ord}}
\nc{\nsu}{\nsubseteq}
\nc{\emp}{\emptyset}

\title{Tropical decomposition of Young's partition lattice}

\author{Vivek Dhand}

\begin{document}

\maketitle
\thispagestyle{empty}

\vspace{-0.2in}

\begin{abstract}
Young's partition lattice $L(m,n)$ consists of unordered partitions having $m$ parts where each part is at most $n$.  Using methods from complex algebraic geometry, R. Stanley proved that this poset is rank-symmetric, unimodal, and strongly Sperner.  Moreover, he conjectured that it has a stronger property called symmetric chain decomposition. Despite many efforts, this conjecture has only been proved for $\min(m,n)\leq 4$.   In this paper, we decompose $L(m,n)$ into level sets for certain tropical polynomials derived from the secant varieties of the rational normal curve in projective space, and we find that the resulting subposets have an elementary raising and lowering algorithm.  As a corollary, we obtain a symmetric chain decomposition for the subposet of $L(m,n)$ consisting of ``sufficiently generic" partitions.
\end{abstract}


\vspace{-0.1in}

\se{Introduction}
Young's partition lattice $L(m,n)$ is defined to be the poset of unordered partitions $\la = (0 \leq \la_1 \leq \dots \leq \la_m \leq n)$ equipped with the following partial order:
\[ \la \leq \mu \iff \la_i \leq \mu_i \;\tn{ for all }\; 1 \leq i \leq m . \]
We can visualize the elements of this poset as Young diagrams which fit in the bottom left-hand corner of an $(m \times n)$ rectangle, ordered by inclusion.  Note that $L(m,n)$ is a  ranked poset, where the rank of $\la$ is given by $\la_1 + \dots + \la_m$.  

This poset appears in different guises in several branches of mathematics. For example, it is isomorphic the poset of Schubert cells in the Grassmannian of $m$-planes in $\bC^{m+n}$.  In the groundbreaking paper \cite{S}, R. Stanley applied the hard Lefschetz theorem to prove that $L(m,n)$ is rank-symmetric, unimodal, and strongly Sperner.  Furthermore, he conjectured that it has a {\em symmetric chain decomposition}, i.e. can be expressed as a disjoint union of rank-symmetric, saturated chains.  Even after many years, this conjecture has only been proved for $\min(m,n)\leq 4$ \cite{L,W}.  In elementary terms, the problem is to find a rule such that:

\begin{enumerate}
\item For each Young diagram in $L(m,n)$, we either do nothing or remove a box so that the result is another Young diagram.
\item Each Young diagram has at most one preimage under this rule.
\item Each terminal Young diagram has complementary rank with the corresponding initial Young diagram.
\end{enumerate}

The main result of this paper involves a simple algorithm which, among other things, yields a symmetric chain decomposition for a large subposet of $L(m,n)$ consisting of ``sufficiently generic" partitions.

\sse*{Example:} $L(2,3)$.
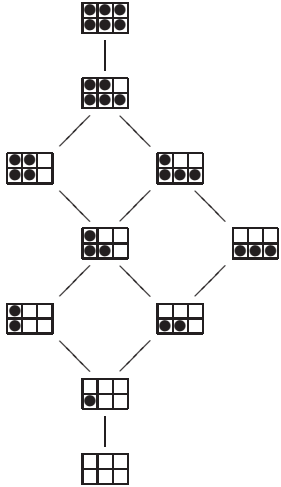
\[\xy
(0,0); (6,0) **\dir{-}; (6,4) **\dir{-}; (0,4) **\dir{-}; (0,0) **\dir{-}; (0,2); (6,2) **\dir{-}; (2,0); (2,4) **\dir{-}; (4,0); (4,4) **\dir{-}; 
(0,10); (6,10) **\dir{-}; (6,14) **\dir{-}; (0,14) **\dir{-}; (0,10) **\dir{-}; (0,12); (6,12) **\dir{-}; (2,10); (2,14) **\dir{-}; (4,10); (4,14) **\dir{-}; 
(10,20); (16,20) **\dir{-}; (16,24) **\dir{-}; (10,24) **\dir{-}; (10,20) **\dir{-}; (10,22); (16,22) **\dir{-}; (12,20); (12,24) **\dir{-}; (14,20); (14,24) **\dir{-}; 
(-10,20); (-4,20) **\dir{-}; (-4,24) **\dir{-}; (-10,24) **\dir{-}; (-10,20) **\dir{-}; (-10,22); (-4,22) **\dir{-}; (-8,20); (-8,24) **\dir{-}; (-6,20); (-6,24) **\dir{-}; 
(0,30); (6,30) **\dir{-}; (6,34) **\dir{-}; (0,34) **\dir{-}; (0,30) **\dir{-}; (0,32); (6,32) **\dir{-}; (2,30); (2,34) **\dir{-}; (4,30); (4,34) **\dir{-}; 
(20,30); (26,30) **\dir{-}; (26,34) **\dir{-}; (20,34) **\dir{-}; (20,30) **\dir{-}; (20,32); (26,32) **\dir{-}; (22,30); (22,34) **\dir{-}; (24,30); (24,34) **\dir{-}; 
(10,40); (16,40) **\dir{-}; (16,44) **\dir{-}; (10,44) **\dir{-}; (10,40) **\dir{-}; (10,42); (16,42) **\dir{-}; (12,40); (12,44) **\dir{-}; (14,40); (14,44) **\dir{-}; 
(-10,40); (-4,40) **\dir{-}; (-4,44) **\dir{-}; (-10,44) **\dir{-}; (-10,40) **\dir{-}; (-10,42); (-4,42) **\dir{-}; (-8,40); (-8,44) **\dir{-}; (-6,40); (-6,44) **\dir{-}; 
(0,50); (6,50) **\dir{-}; (6,54) **\dir{-}; (0,54) **\dir{-}; (0,50) **\dir{-}; (0,52); (6,52) **\dir{-}; (2,50); (2,54) **\dir{-}; (4,50); (4,54) **\dir{-}; 
(0,60); (6,60) **\dir{-}; (6,64) **\dir{-}; (0,64) **\dir{-}; (0,60) **\dir{-}; (0,62); (6,62) **\dir{-}; (2,60); (2,64) **\dir{-}; (4,60); (4,64) **\dir{-}; 
(1,61) *+{\bullet}; (3,61) *+{\bullet}; (5,61) *+{\bullet}; (1,63) *+{\bullet}; (3,63) *+{\bullet}; (5,63) *+{\bullet};
(1,51) *+{\bullet}; (3,51) *+{\bullet}; (5,51) *+{\bullet}; (1,53) *+{\bullet}; (3,53) *+{\bullet}; 
(11,41) *+{\bullet}; (13,41) *+{\bullet}; (15,41) *+{\bullet}; (11,43) *+{\bullet}; 
(-9,41) *+{\bullet}; (-7,41) *+{\bullet}; (-9,43) *+{\bullet}; (-7,43) *+{\bullet}; 
(1,31) *+{\bullet}; (3,31) *+{\bullet}; (1,33) *+{\bullet}; 
(21,31) *+{\bullet}; (23,31) *+{\bullet}; (25,31) *+{\bullet}; 
(-9,21) *+{\bullet}; (-9,23) *+{\bullet}; 
(11,21) *+{\bullet}; (13,21) *+{\bullet};
(1,11) *+{\bullet}; 
(3,5); (3,9) **\dir{-};
(1,15); (-3,19) **\dir{-}; (5,15); (9,19) **\dir{-};
(-3,25); (1,29) **\dir{-}; (9,25); (5,29) **\dir{-}; (15,25); (19,29) **\dir{-};
(1,35); (-3,39) **\dir{-}; (5,35); (9,39) **\dir{-}; (19,35); (15,39) **\dir{-};
(-3,45); (1,49) **\dir{-}; (9,45); (5,49) **\dir{-};
(3,55); (3,59) **\dir{-};
\endxy\]


\sse{The monomial model of Young's partition lattice} 

Note that a partition $\la \in L(m,n)$ is uniquely determined by the $(n+1)$-tuple $(a_0, \dots, a_n)$ where $a_i$ is equal to the number of times that $i$ appears in $(\la_1, \dots, \la_m)$.  A covering relation in $L(m,n)$ is given by adding a box to an acceptable row of a Young diagram, which sends:
\[ (a_0, \dots, a_n) \mapsto (a_0, \dots,a_{i-1}-1,a_i+1, \dots, a_n). \]
In this way, we obtain an edge coloring of the Hasse diagram of $L(m,n)$, using $n$ colors, where the $i$-th color corresponds to the above operation.   

Let $A_n = \bC[z_0, \dots, z_n]$. By mapping $(a_0, \dots, a_n)$ to the monomial $z_0^{a_0} \dots z_n^{a_n}$, we see that $L(m,n)$ is isomorphic to the poset $A_n(m)$ of monomials of degree $m$ in $A_n$.  In terms of partitions, the isomorphism is given by:
\[  (\la_1, \dots, \la_m) \mapsto z_{\la_1} z_{\la_2} \dots z_{\la_m}. \]
The $i$-th color operation corresponds to changing a single $z_{i-1}$ to $z_{i}$, and this partial order is induced by the standard action of $\slc$ on the irreducible representation $\bC \lan z_0, \dots, z_n \ran$.
The advantage of using this model for Young's partition lattice is that we can simultaneously deal with $A_n(m)$ for all $m \geq 0$ and exploit the full structure of $A_n$ as a commutative graded $\bC$-algebra with $\slc$-action.  

\sse{Combinatorial secant ideals of the rational normal curve}

Since finding a symmetric chain decomposition of $L(m,n)$ is so difficult, we might first look for a coarser decomposition into centered $\slc$-subposets.  From the point of view of commutative algebra, a natural way to decompose the monomial basis in $A_n$ is to choose a monomial ideal $I \su A_n$ and form the graded algebra:
\[ \gr_{I}(A_n) = \bigoplus_{j \geq 0} I^j/I^{j+1}.  \]
In fact, we can obtain even finer decompositions by starting with a set of monomial ideals $I_1, \dots, I_t \su A_n$ and looking at the multi-graded algebra:
\[ \gr_{I_1} \dots \gr_{I_t} A_n. \] 
In order to take full advantage of the underlying symmetry, we should start with a set of $\slc$-invariant ideals and perform a Gr\"obner degeneration to obtain monomial ideals.  

From the point of view of algebraic geometry, the nicest possible $\slc$-invariant subvarieties of projective space are the {\em rational normal curve} $\cC_n \su \bP^n$ and its {\em secant varieties} $\cC_n^{\{r\}}$ for $0 \leq r \leq k = \lfloor n/2 \rfloor$.  It is well known that the ideal of $\cC_n^{\{r\}}$ is generated by maximal minors of the $(r + 1) \times (n+1-r)$ Hankel matrix:
\[ \left[ \begin{array}{cccc} z_0 & z_1 & \dots & z_{n-r} \\ z_1 & z_2 & \dots & z_{n+1-r} \\ \vdots & \vdots & & \vdots \\ z_r & z_{r+1} &\dots & z_n \end{array} \right]. \] 
By considering the initial ideals of $\cC_n^{\{r\}}$ (with respect to any diagonal term order)  we obtain a set of squarefree monomial ideals $I_{n,r} \su A_n$, for $0 \leq r \leq k$.   

We can now equip the monomial basis in $A_n$ with extra gradings according to the order of vanishing on this set of monomial ideals.  Given a monomial $\mu \in A_n$ and $0 \leq r \leq k$, let $\deg_r(\mu)$ be the number of minimal generators of $I_{n,r}$ appearing in the image of $\mu$ under the natural map:
\[  A_n \to \gr_{I_{n,1}} \dots \gr_{I_{n,k}} A_n . \]
Given integers $d_0, \dots, d_k \geq 0$, we define:
\[ Q_n(d_0, \dots, d_k) = \{ \mu \in A_n \mid \deg_r(\mu) = d_r \;\tn{ for each }\; 0 \leq r \leq k  \}. \]
Set-theoretically, this decomposition is implicit in Conca's ``canonical decomposition" algorithm for factoring monomials \cite{Co}.  However, our description (in terms of tropical polynomials) does not depend on a particular choice of factorization and therefore allows a finer investigation of the poset structure.

It turns out that $I_{n,1}$ is the Stanley-Reisner ideal of the path graph with vertices $\{0,\dots, n\}$, and $I_{n,r}$ is the $r$-th combinatorial secant ideal of $I_{n,1}$.  Using this fact, we find an explicit formula for the irredundant irreducible decomposition of each $I_{n,r}$, which leads us to the following tropical polynomials:
\[ f_{n,r}(a_0, \dots, a_n) =  \min_{0\leq \la_0 \leq \dots \leq \la_{n-2r} \leq r}\sum_{j = 0}^{n-2r} a_{2\la_j + j}. \]
For any monomial  $\mu = z_0^{a_0} \dots z_n^{a_n}$, we have:
\[ \mu \in I_{n,r}^{(s)} \iff f_{n,r}(\mu) \geq s, \]
where $I_{n,r}^{(s)}$ denotes the $s$-th symbolic power and $f_{n,r}(\mu)$ denotes $f_{n,r}(a_0, \dots, a_n)$.  Remarkably, grading by the symbolic powers give us the same decomposition of monomials as the ordinary powers:
\[ Q_n(d_0, \dots, d_k) = \{ \mu \in A_n \mid  f_{n,r}(\mu) = \sum_{j = r}^{k} ( j + 1 - r ) d_j   \;\tn{ for each }\; 0 \leq r \leq k  \}. \]
With this explicit description in hand, we find that $Q_{n-2}(d_1, \dots, d_k)$ embeds into $Q_n(d_0, \dots, d_k)$ in two opposite ways, and there exists an elementary algorithm for producing coverings by saturated chains running between these two extremes.

\sse{The raising and lowering algorithm} We can associate a sequence of colors to any saturated chain in $A_n(m)$ by reading the covering relations along the chain from highest to lowest weight.  A saturated chain is said to be {\em monotonic} if its color sequence is monotonically increasing.  

Given a monomial $\mu = z_0^{a_0} \dots z_n^{a_n} \in Q_n(d_0, \dots, d_k)$, the raising and lowering algorithm produces several monotonic saturated chains in $Q_n(d_0, \dots, d_k)$ passing through $\mu$.  A pair of adjacent entries $(a_i,a_{i+1})$ is called a {\em maximal pair} if:
\[ a_i + a_{i+1} = \underset{0 \leq j \leq n-1}{\tn{max}} (a_j + a_{j+1}). \]
The steps of the algorithm can be broken down into elementary moves of three types:

({\bf cover}) Choose a maximal pair $(a_i,a_{i+1})$  from $(a_0, \dots, a_n)$.

({\bf move right}) Apply the $(i+1)$-st color to get $(a_i-1,a_{i+1}+1)$.

({\bf move left}) Apply the inverse of the $(i+1)$-st color to get $(a_i+1,a_{i+1}-1)$. 

\sse*{The algorithm (right-moving version):}
\begin{enumerate}
\item Cover a maximal pair $(a_i, a_{i+1})$ for some $0 \leq i \leq n-1$. If $i = n-1$, move right $a_{n-1}$ times and end the chain.
\item Move right $(a_i - a_{i+2})$ times, and then go back to step (1) and cover $(a_{i+1}, a_{i+2})$. 
\end{enumerate}

\sse*{The algorithm (left-moving version):}
\begin{enumerate}
\item Cover a maximal pair $(a_{i-1}, a_i)$ for some $1 \leq i \leq n$. If $i = 1$, move left $a_1$ times and end the chain.
\item Move left $(a_i - a_{i-2})$ times, and then go back to step (1) and cover $(a_{i-2}, a_{i-1})$.
\end{enumerate}

While the algorithm potentially outputs several saturated chains involving the monomial $\mu$, there are two chains in particular that are important for us.  If we start at the leftmost coverable pair, we can combine the two versions of the algorithm to get a single monotonic saturated chain. The highest and lowest weight elements of this chain will respectively satisfy the following equations:
\[ a_0 =  \underset{0 \leq j \leq n-1}{\tn{max}} (a_j + a_{j+1}) \;\;,\;\; a_n =  \underset{0 \leq j \leq n-1}{\tn{max}} (a_j + a_{j+1}). \]
Alternately, we could start with the rightmost coverable pair and get another monotonic saturated chain whose respective endpoints satisfy the above equations. We will call the two monotonic saturated chains obtained in this way the left and right {\em transversal chains} of $\mu$.  

\sse*{Theorem A} Let $n, d_0, \dots, d_k \geq 0$. The transversal chains provide a vertex covering of $Q_n(d_0, \dots, d_k)$ by monotonic saturated chains.  If $d_0 > 0$, then this covering is actually a decomposition.

Transversal chains are not necessarily rank-symmetric.  Nevertheless, when dealing with sufficiently generic monomials, we can stitch them together to get a symmetric chain decomposition.

\sse*{Theorem B} If $Q_{n-2}(d_1, \dots, d_k)$ has a symmetric chain decomposition and
\[ 1+ 2 d_0 \geq  \displaystyle \sum_{i = 1}^{k} d_i (n-2i)i, \]
then $Q_n(d_0,\dots,d_k)$ has a symmetric chain decomposition.


In other words, we can decompose Young's partition lattice into a ``generic" part and a ``singular" part:
\[ L(m,n) = L(m,n)^{gen} \sqcup L(m,n)^{sing} \]
where $L(m,n)^{gen}$ is a symmetric chain order and $L(m,n)^{sing}$ has a covering by monotonic saturated chains. Note that by a decomposition of a poset we mean a set-theoretic decomposition where each factor is equipped with the induced partial order.



Let us outline the contents of the paper.  In section 2, we decompose $L(m,n) \simeq A_n(m)$ according to the order of vanishing on the monomial ideals $I_{n,r}$, for $0 \leq r \leq \lfloor n/2\rfloor$, and we give a simpler description in terms of the level sets of certain tropical polynomials.  The poset structure aside, the results of this section are essentially contained in \cite{Co}.  In section 3, we use the tropical description to prove that each level set in $A_n$ contains two extremal embeddings of a smaller level set in $A_{n-2}$.  In section 4, we describe the raising and lowering algorithm and prove that it yields a covering of the level sets by monotonic saturated chains running between the extremal subposets.  Indeed, most of the time this covering is a decomposition, and this fact allows us to inductively stitch these chains together to obtain a symmetric chain decomposition for the generic part of $L(m,n)$.
 

\sse*{Acknowledgements}  Thanks are due to Peter Magyar for telling me about Stanley's conjecture, and for many enlightening conversations on a host of mathematical topics.  I would also like to thank the Department of Mathematics at Michigan State University for their hospitality while this project was under way.

\se{Combinatorial secant ideals of the rational normal curve}

Let $A_n = \bC[z_0,\dots, z_n]$ and let $A_n(m)$ denote the set of monomials of degree $m$ in $A_n$.  Define a partial order on $A_n(m)$ as follows:
\[ z_0^{a_0} \dots z_n^{a_n} \leq z_0^{b_0} \dots z_n^{b_n} \iff \sum_{i = j}^n a_i \leq \sum_{i = j}^n b_i \tn{ for } 1 \leq j \leq n . \]
This partial order is induced by the following action of $\slc = \bC \lan H,E, F\ran$ on the irreducible representation $\bC \lan z_0, \dots, z_n \ran$:
\[ H(z_i) = n - 2i \quad\quad  E(z_i) = i z_{i-1} \quad\quad  F(z_i) = (n-i) z_{i+1}. \]
It is not hard to see that $L(m,n)$ is isomorphic to $A_n(m)$ by the map:
\[  (\la_1, \dots, \la_m) \mapsto z_{\la_1} z_{\la_2} \dots z_{\la_m}. \]
Note that $L(m,n) \simeq L(n,m)$ simply by taking the conjugate partition, which corresponds to the well-known duality $A_n(m) \simeq A_m(n)$ for choosing multisets from a finite set. Also note that the weight of a monomial $\mu = z_0^{a_0} \dots z_n^{a_n}$ is given by:
\[ \tn{wt}(\mu) = \sum_{i = 0}^n a_i (n-2i), \]
while its rank is given by:
\[ \tn{rk}(\mu) = \sum_{i = 0}^n i a_i, \]
and therefore we obtain the simple relation:
\[  \tn{wt}(\mu) =  m n  - 2 \tn{rk}(\mu). \]
Since we will be working exclusively with centered subposets of $A_n(m)$, it will be more convenient for us to refer to the weights of monomials instead of their ranks.

Recall that the ideal of the rational normal curve $\cC_n \su \bP^n$ is generated by the set of maximal minors of the $(2 \times n)$-Hankel matrix:
\[  H_{n,1} = \left[ \begin{array}{cccc} z_0 & z_1 & \dots & z_{n-1} \\ z_1 & z_2 & \dots & z_{n} \end{array} \right] . \]
Furthermore, the ideal of the $r$-th secant variety of $\cC_n$ is generated by the set of maximal minors of the $(r+1) \times(n-r+1)$ Hankel matrix:
\[  H_{n,r} = \left[ \begin{array}{cccc} z_0 & z_1 & \dots & z_r \\ z_1 & z_2 & \dots & z_{r+1} \\ \vdots & \vdots && \vdots \\ z_{n-r} & z_{n-r+1}& \dots & z_n \end{array} \right] \]
where $1 \leq r \leq \lfloor n/2 \rfloor$. Let $I_{n,r}$ denote the initial ideal (with respect to any diagonal term order) of the ideal of the $r$-th secant variety of $\cC_n$.  The minimal generators of $I_{n,r}$ are the initial monomials of Hankel determinants:
\[ I_{n,r} = ( \{ z_{i_0} \dots z_{i_{r}} \mid  i_j +1 < i_{j+1} \tn{ for } 0 \leq j \leq r-1  \} ) . \]
For ease of notation, we define $I_n = I_{n,1}$ and $ I_{n,0} = \fm = (z_0, \dots, z_n) $. 

Recall that for an ideal $I$ in a ring $A$ and an element $a \in A$ such that:
\[ a \notin \bigcap_{j \geq 0} I^j, \]
the {\em order of vanishing of} $a$ {\em on} $I$ is defined as:
\[ \ord_{I}(a) = \max \{ j \mid a \in I^j  \} . \]  
The associated graded ring
\[ \gr_I A = \bigoplus_{k \geq 0} I^k/I^{k+1} \]
has a {\em symbol map}
\[ \si: A \to \gr_I A \]
where:
\[ \si(a) = a \;(\tn{mod } I^{d + 1}) \;\;\tn{if}\;\;\ord_I(a) = d, \]
\[ \si(a) = 0 \;\;\tn{if}\;\; a \in \bigcap_{j \geq 0} I^j . \]
Also, note that if $J \su A$ is another ideal, we can form the multi-graded algebra:
\[ \gr_J \gr_I A :=  \gr_{\wt{J}} \; \gr_I A, \]
where $\wt{J}$ denotes the ideal generated by $\si(J) \su \gr_I(A)$.

In our situation, we will consider the consider the composition of several symbol maps:
\[ A_n \to \gr_{I_{n,k}} A_n \to  \gr_{I_{n,k-1}} \gr_{I_{n,k}} A_n \to \dots \to \gr_{I_{n,1}} \dots  \gr_{I_{n,k}} A_n. \]
Let us describe the steps involved in calculating this map. Given a monomial $\mu \in A_n(m)$, let $\deg_k(\mu) =\ord_{I_{n,k}}(\mu)$.   From $\mu$, we factor out a monomial $\mu_k$ which is a product of $\deg_k(\mu)$ minimal generators of $I_{n,k}$.  Now let $\deg_{k-1}(\mu) = \ord_{I_{n,k-1}}(\mu/\mu_k)$ and factor out from $\mu/\mu_k$ a monomial $\mu_{k-1}$ which is a product of $\deg_{k-1}(\mu)$ minimal generators of $I_{n,k-1}$.  Continuing in this way, we obtain a {\em maximal factorization} of $\mu$:
\[ \mu = \mu_0 \mu_1  \dots \mu_k. \]
While this factorization is not necessarily unique, what matters is that $\deg_r(\mu)$ is uniquely determined by $\mu$, for all $0 \leq r \leq k$.  

We define the following decomposition of the monomial basis in $A_n$:
\[  Q_n(d_0, \dots, d_k) = \{  \mu \in A_n \mid \deg_r(\mu) = d_r  \tn{ for all } 0 \leq r \leq k \}. \]
Note that $Q_n(d_0, \dots, d_k) \su A_n(m)$, where:
\[ m = \sum_{r = 0}^k d_r(r+1) . \]
We can visualize a maximal factorization as a tableau with $\deg_r(\mu)$ rows of length $r+1$, where the entries in each row form a $<_1$-chain, i.e. a sequence of non-negative integers $(i_0, \dots, i_r)$ such that  $i_j +1 < i_{j+1} $ for $ 0 \leq j \leq r-1$.  We can rearrange the entries in this tableau so that each entry is as small as possible as we read from left to right and from the longest to the shortest row.  In this way we obtain the same tableau as Conca's ``canonical decomposition" of monomials \cite{Co}.  
While our decomposition is set-theoretically equivalent to Conca's decomposition, there are a few subtle but important differences.  Conca's construction yields one particular maximal factorization for each monomial in terms of minimal generators of the ideals $I_{n,r}$.  This choice obscures the poset structure, which is our main object of interest.  We will give a description of $Q_n(d_0, \dots, d_k)$ which treats all maximal factorizations equally and illuminates the structural relationships between these posets.

First, we need an explicit description of the (unique) irredundant irreducible decomposition of each $I_{n,r}$.  Our description will be in terms of the associated simplicial complex.  Let $\De_n$ denote the set of subsets of $\{0, \dots, n\}$.  For any $F \in \De_n$, we define
\[ z^F = \prod_{i \in F} z_i \in A_n \quad\tn{and}\quad \fm^F = \sum_{i \in F} (z_i) \su A_n .\]
Let $\Ga \su \De_n$ be an abstract simplicial complex, so $S \su T \in \Ga \implies S \in \Ga$.   The Stanley-Reisner ideal of $\Ga$ is defined as \cite{MS}:
\[ I_\Ga = \lan z^F \mid F \notin \Ga  \ran  = \bigcap_{\ov{F} \in \Ga} \fm^F. \]
For our problem, the relevant simplicial complex  is the path graph with $n+1$ vertices:
\[ \Ga_n = \{\emp, \{0\}, \dots, \{n\}, \{0,1\}, \{1,2\}, \dots, \{n-1,n\} \} . \]
It is easy to check that the Stanley-Reisner ideal of $\Ga_n$ is equal to $I_n$.  Indeed, $I_n \su I_{\Ga_n}$ because $I_n$ is generated by the quadratic monomials corresponding to the edges in the complement of the path graph.  On the other hand, any $F \in \De_n$ such that $|F| \geq 3$ must contain some $i$ and $j$ such that $i+1 < j$, so we see that $I_{\Ga_n} \su I_n$.





\sse{Remark}
We know from Section 6.1 of \cite{Su} that $I_{n,r}$ is equal to the $r$-th combinatorial secant ideal of $I_n$.  Let $\Ga_{n,r}$ denote the simplicial complex of $r$-fold unions of simplices from $\Ga_n$:
\[ \Ga_{n,r} = \{ F_1 \cup \dots \cup F_r \mid F_i \in \Ga_n \}. \]
By Remark 2.9 of \cite{StSu}, it follows that $I_{n,r}$ is the Stanley-Reisner ideal of $\Ga_{n,r}$.  Furthermore, it is not hard to see that each facet of $\Ga_{n,r}$ is equal to the disjoint union of $r$ edges in $\Ga_n$.

\sse{Proposition} The irredundant irreducible decomposition of $I_{n,r}$ is given by:
\[ I_{n,r} = \bigcap_{0\leq \la_0 \leq \dots \leq \la_{n-2r} \leq r} \fm^{\{2\la_0, 2\la_1 + 1 , \dots, 2\la_{n-2r} + n-2r  \}} . \]

\begin{proof} 
From the definition of Stanley-Reisner ideal, we have:
\[ I_{n,r} = \bigcap_{\ov{F} \in \Ga_{n,r}} \fm^{F}. \]
It suffices to take the intersection over the set of complements of facets of $\Ga_{n,r}$:
\[  \{ \ov{T_1 \sqcup \dots \sqcup T_r} \mid T_1, \dots, T_r \tn{ are disjoint edges of } \Ga_n \}.  \]
If we list the elements of $\ov{T_1 \sqcup \dots \sqcup T_r}$ in increasing order, we see that each term will be sandwiched between sequences of consecutive edges chosen from $\Ga_n$.  Therefore, the set of $F \in \De_n$ such that $\ov{F} \in \Ga_{n,r}$ is equal to:
\[  \{ \{ 2\la_0, 2\la_1 + 1 , \dots, 2\la_{n-2r} + n-2r \} \mid 0 \leq \la_0 \leq \dots \leq \la_{n-2r} \leq r \} \]
and the result follows (cf. \cite{Co}, Lemma 3.5).
\end{proof}

We are now ready to describe the symbolic powers of our monomial ideals.   If $I$ is a radical ideal in a polynomial ring over an algebraically closed field, then the $s$-th {\em symbolic power of} $I$ is:
\[  I^{(s)} = \bigcap_{M \in \cM_I} M^s \]
where $\cM_I$ denotes the set of all maximal ideals containing $I$ \cite{Su}.

\sse{Proposition} For each $0 \leq r \leq \lfloor n/2 \rfloor$, the tropical polynomial:
\[ f_{n,r}(a_0, \dots, a_n) =  \min_{0\leq \la_0 \leq \dots \leq \la_{n-2r} \leq r}\sum_{j = 0}^{n-2r} a_{2\la_j + j} \]
satisfies the following property:
\[ z_0^{a_0} \dots z_n^{a_n} \in I_{n,r}^{(s)} \iff f_{n,r}(a_0, \dots, a_n) \geq s . \]
\begin{proof}
In fact, there is a general relationship between symbolic powers and tropical polynomials. Any squarefree monomial ideal $I$ with the irredundant irreducible decomposition
\[ I = \bigcap_{F} \fm^F \]
has the following symbolic powers:
\[ I^{(s)} = \bigcap_{F} (\fm^{F})^s . \]
Therefore,
\[ z_0^{a_0} \dots z_n^{a_n} \in I^{(s)}  \iff  z_0^{a_0} \dots z_n^{a_n} \in (\fm^F)^s \tn{ for all } F .\]
Now,
\[ (\fm^F)^s = \left( \left\{ \prod_{i \in F} z_i^{b_i} \;\bigg{|}\; \sum_{i\in F} b_i = s \right\} \right) \]
and so:
\[  z_0^{a_0} \dots z_n^{a_n} \in (\fm^F)^s  \iff  \sum_{i\in F} a_i \geq s . \]
In our particular example, we obtain:
\[ z_0^{a_0} \dots z_n^{a_n} \in I_{n,r}^{(s)} \iff \sum_{j = 0}^{2n-r} a_{2\la_j + j} \geq s \tn{ for all } 0\leq \la_0 \leq \dots \leq \la_{n-2r} \leq r , \]
which is equivalent to $f_{n,r}(a_0, \dots, a_n) \geq s.$
\end{proof}

\sse{Example} 




%



%





If $n = 5$, then the ideals:
\[ I_{5,1} = (z_0 z_2,z_0 z_3,z_0 z_4,z_0 z_5,z_1 z_3,z_1 z_4,z_1 z_5,z_2 z_4,z_2 z_5,z_3 z_5) \] 
\[ I_{5,2} = (z_0z_2z_4,z_0z_2z_5,z_0z_3z_5,z_1z_3z_5)\]
have the following irredundant irreducible decompositions:
\[ I_{5,1} = (z_0,z_1,z_2,z_3) \cap (z_0,z_1,z_2,z_5) \cap (z_0,z_1,z_4,z_5) \cap (z_0,z_3,z_4,z_5) \cap (z_2,z_3,z_4,z_5) \]
\[ I_{5,2} = (z_0,z_1)\cap(z_0,z_3)\cap(z_0,z_5)\cap(z_2,z_3)\cap(z_2,z_5)\cap(z_4,z_5) \]
and the corresponding tropical polynomials are: 
\[ f_{5,1} \!=\! \min(a_0+a_1+a_2+a_3, a_0 + a_1 + a_2 + a_5, a_0 + a_1 + a_4 + a_5, a_0 + a_3 + a_4 + a_5, a_2 + a_3 + a_4 + a_5 ) \]
\[ f_{5,2} = \min(a_0 + a_1, a_0 + a_3, a_0 + a_5,a_2 + a_3, a_2 + a_5, a_4 + a_5) . \]

\sse{Remark}  For a monomial $\mu = z_0^{a_0} \dots z_n^{a_n}$, note that $f_{n,r}$ minimizes a sum over the complements of facets of $\Ga_{n,r}$:
\[ f_{n,r}(\mu) = \underset{F \in \Ga_{n,r}}{\min} \sum_{i \notin F} a_i. \]
Since $\deg(\mu) = a_0 + \dots + a_n$, we can rewrite this as:
\[ f_{n,r}(\mu) = \underset{F \in \Ga_{n,r}}{\min} (\deg(\mu) - \sum_{i \in F} a_i) = \deg(\mu) - \underset{{F \in \Ga_{n,r}}}{\max} \sum_{i \in F} a_i.  \]
Given a facet $F \in \Ga_{n,r}$, we will refer to the following sum: 
\[ \al(F,\mu) = \sum_{i \in F} a_i \] as the {\em amount of} $\mu$ {\em covered by} $F$.
With this terminology, in order to calculate $f_{n,r}(\mu)$ we need to use $r$ disjoint edges in $\Ga_n$ to cover as much of $\mu$ as possible.

\sse{Lemma} For any monomial $\mu \in A_n$ and $0 \leq r \leq \lfloor n/2 \rfloor$, we have:
\[  f_{n,r}(\mu) = \sum_{j = r}^{\lfloor n/2 \rfloor} ( j + 1 - r ) \deg_j(\mu) \;\;\tn{and}\;\; \deg_r(\mu) = f_{n,r}(\mu) - 2 f_{n,r+1}(\mu) + f_{n,r+2} ,\] 
where $f_{n,j} = 0$ for $j> \lfloor n/2 \rfloor$.

\begin{proof}
Recall that a maximal factorization of $\mu$ can be expressed as a tableau with $\deg_r(\mu)$ rows containing $<_1$-chains of size $(r+1)$.  If we rearrange the entries in this tableau so that each entry is as small as possible as we read from left to right and from the longest to the shortest row, we obtain the same result as factoring out the generators of $I_{n,r}$ in lexicographic order, as $r$ runs from $k$ to $0$ (i.e. Conca's ``canonical decomposition"). The key property of this ordering is that, for any entry $a$, either $a$ or $a-1$ must appear in each preceding row.  If not, we could insert $a$ earlier in the tableau while preserving the $<_1$-chains in the rows, which would contradict the minimality of the ordering \cite{Co}.  

Now consider the smallest entry $a$ in the shortest row of the tableau.  Since $a$ or $a-1$ must appear in each preceding row, we see that covering $\{a-1,a\}$ with an edge from $\Ga_n$ will maximize the amount covered.  Moreover, the rows are $<_1$-chains, so we can cover exactly one entry from each row.  In other words, the number of boxes in the first column of the tableau is equal to:
\[ \underset{{F \in \Ga_{n}}}{\max} \sum_{i \in F} a_i . \]
Therefore, $f_{n,1}(\mu)$ is equal to the number of boxes that do not lie in the first column.

Now we remove the boxes with entries $a$ or $a-1$ from the tableau and repeat this argument for what remains.  It follows that $f_{n,r}(\mu)$ is equal to the number of boxes of the tableau that do not lie in the first $r$ columns.  Subtracting $r$ boxes from each row and adding up what remains, we find that:
\[ f_{n,r}(\mu) = \sum_{j = r}^{k} ( j + 1 - r ) \deg_j(\mu). \]
Moreover, $\deg_r(\mu)$ is equal to the difference between the number of boxes in the $(r+1)$-th column and the $(r+2)$-th column. Therefore:
\begin{eqnarray*}
 \deg_r(\mu) &=& (f_{n,r}(\mu) - f_{n,r+1}(\mu)) - (f_{n,r+1}(\mu) - f_{n,r+2}(\mu))\\
 &=& f_{n,r}(\mu) - 2f_{n,r+1}(\mu) + f_{n,r+2}(\mu)
\end{eqnarray*}
as desired.
\end{proof}


\sse{Example} Consider the monomial $\mu = z_0^4 z_1^3 z_2^2 z_3 z_4 z_5^4 \in A_{5}(15).$  Conca's decomposition for $\mu$ yields the tableau:
\[ \xy (0,0); (0,35) **\dir{-}; (0,0); (5,0) **\dir{-}; (5,35) **\dir{-}; (0,5); (5,5) **\dir{-}; (0,10); (10,10) **\dir{-}; (10,35) **\dir{-}; (0,15); (10,15) **\dir{-}; (0,20); (15,20) **\dir{-}; (15,35) **\dir{-}; (0,25); (15,25) **\dir{-};
(0,30); (15,30) **\dir{-};
(0,35); (15,35) **\dir{-};
(2.5,2.5) *+{1};
(2.5,7.5) *+{1};
(2.5,12.5) *+{1};
(2.5,17.5) *+{0};
(2.5,22.5) *+{0};
(2.5,27.5) *+{0};
(2.5,32.5) *+{0};
(7.5,12.5) *+{5};
(7.5,17.5) *+{5};
(7.5,22.5) *+{3};
(7.5,27.5) *+{2};
(7.5,32.5) *+{2};
(12.5,22.5) *+{5};
(12.5,27.5) *+{5};
(12.5,32.5) *+{4};
 \endxy  \]
Note that $f_{5,0}(\mu) = \deg(\mu) = 15$.  First we cover $\{0,1\}$, which removes $7$ boxes, so $f_{5,1}(\mu) = 15 - 7 = 8$.  Next, we cover $\{4,5\}$, which removes $5$ boxes, so $f_{5,2}(\mu) = 8 - 5 = 3$.   Counting the number of rows of each type, we see that:
\[ \deg_0(\mu) = 2 = 15 - 2(8) + 3 \;\;,\;\; \deg_1(\mu) = 2 = 8 - 2(3) \;\;,\;\;  \deg_2(\mu) = 3, \]
as expected.

\se{Structure of the level sets}


We have the following decomposition of $A_n(m)$:
\[ Q_n(d_0, \dots, d_k) = \{ \mu \in A_n \mid  f_{n,r}(\mu) = \sum_{j = r}^{k} ( j + 1 - r ) d_j   \;\tn{ for each }\; 0 \leq r \leq k  \} \] 
where
\[ m = \sum_{j = 0}^k d_j(j+1) . \]
With this description, we find that these posets have several interesting structures. First, since $\Ga_n$ has an involution which sends $i$ to $n-i$, we see that $I_{n,r}$ and $f_{n,r}$ are invariant under the following automorphism:
\[ \tau_n: A_n \to A_n  \]
\[ z_i \mapsto z_{n-i} . \]
It follows that $\tau_n$ restricts to a rank-flipping involution of $Q_n(d_0, \dots, d_k)$.  In particular, it is a rank-symmetric, centered subposet of $A_n(m)$.  

Also, note that $I_{n,r}$ has a unique minimal generator of highest (resp. lowest) weight, namely:
\[ \mu_{n,r} = z_0 z_2 \dots z_{2r} \;\; (\tn{resp.}\;\; \nu_{n,r} = \tau_n(\mu_{n,r}) = z_n z_{n-2} \dots z_{n-2r}). \]
It follows that $Q_n(d_0, \dots, d_k)$ has a unique monomial of highest (resp. lowest) weight, namely:
\[  \mu_{n}(d_0, \dots, d_k) = \prod_{j = 0}^k \mu_{n,j}^{d_j} \;\; (\tn{resp.} \;\; \nu_{n}(d_0, \dots, d_k) = \prod_{j = 0}^k \nu_{n,j}^{d_j}).  \]
We can now state the embedding property of the $Q$-posets.

\sse{Lemma} \label{Embedding} Let $d_0, \dots, d_k \geq 0$ and let $d = d_0 + \dots + d_k$.  
Then multiplication by $z_{n}^{d}$ induces an embedding of posets:
\[ z_n^d: Q_{n-2}(d_1,\dots,d_k) \to Q_n(d_0, \dots, d_k) \]
which sends $\nu_{n-2}(d_1, \dots, d_k)$ to $\nu_{n}(d_0, \dots, d_k)$.

\begin{proof}
Since $z_n \nu_{n-2,r} = \nu_{n,r}$, we see that:
\[ z_n^d \nu_{n-2}(d_1, \dots, d_k) = z_n^{d_0} \prod_{j = 1}^{k} z_{n}^{d_j} \nu_{n-2,j}^{d_j} = \prod_{j = 0}^k \nu_{n,j}^{d_j} = \nu_n(d_0, \dots, d_k). \]
Let $\mu \in Q_{n-2}(d_1, \dots, d_k)$.  Recall that: 
\[ (f_{n,0} - f_{n,1})(a_0, \dots, a_n) = \max_{0 \leq j \leq n-1} (a_j + a_{j+1}). \]
In the monomial $z_n^d \mu$, the sum of the exponents of $z_{n-1}$ and $z_{n}$ is equal to $d$, and:
\[ d \geq d_1 + \dots + d_k = f_{n-2,0}(\mu) - f_{n-2,1}(\mu). \]  
In other words, the edge $\{n-1,n\}$ will cover at least as much of $z_n^d \mu$ as any other edge in $\Ga_n$.

Now we claim that, for each $1 \leq r \leq k$:
\[ f_{n,r}(z_n^d \mu) = f_{n-2,r-1}(\mu). \]
Recall that the calculation of $f_{n,r}$ involves all $r$-fold disjoint unions of edges in $\Ga_n$.  We can split up this calculation into two cases: the $n$-th vertex is either covered or uncovered.   If the $n$-th vertex is covered, then so is the $(n-1)$-th vertex, which leaves $r-1$ edges for the vertices $0, \dots, n-2$.  In other words, for $z_n^d \mu$, we have
\[ \max_{n \in F \in \Ga_{n,r}} \sum_{i \in F} a_i =  d + \max_{F \in \Ga_{n-2,r-1}} \sum_{i \in F} a_i = d + \deg(\mu) - f_{n-2,r-1}(\mu). \]
On the other hand, if the $n$-th vertex is uncovered, choose $r-1$ pairwise disjoint edges from $\Ga_{n-2}$, and then arbitrarily choose the $r$-th edge from $\Ga_{n-1}$.  Since the exponent of $z_{n-1}$ in $z_n^d \mu$ is zero, we get the following inequalities:
\begin{eqnarray*} 
\max_{n \notin F \in \Ga_{n,r}} \sum_{i \in F} a_i &\leq& \max_{0 \leq j \leq n-2} (a_j + a_{j+1}) + \max_{F \in \Ga_{n-2,r-1}} \sum_{i \in F} a_i \\
& \leq & d + \deg(\mu) - f_{n-2,r-1}(\mu). 
\end{eqnarray*}
Therefore:
\begin{eqnarray*} 
f_{n,r}(z_n^d \mu) &=& \deg(z_n^d \mu) - \max_{F \in \Ga_{n,r}} \sum_{i \in F} a_i \\
&=& d + \deg(\mu) - \max\left( \max_{n \in F \in \Ga_{n,r}} \sum_{i \in F} a_i , \max_{n \notin F \in \Ga_{n,r}} \sum_{i \in F} a_i \right) \\
&=& d + \deg(\mu) - d - \deg(\mu) + f_{n-2,r-1}(\mu) \\
&=& f_{n-2,r-1}(\mu).
\end{eqnarray*}
It follows that, for $1\leq r \leq k$:
\begin{eqnarray*}
\deg_r(z_n^d \mu) &=& (f_{n,r} - 2 f_{n,r+1} + f_{n,r+2})(z_n^d \mu) \\ &=& (f_{n,r-1} - 2 f_{n,r} + f_{n,r+1})(\mu) \\ &=& \deg_{r-1}(\mu).
\end{eqnarray*}
The only difference occurs at $r = 0$:
\begin{eqnarray*}
\deg_0(z_n^d \mu) &=& (f_{n,0} - 2 f_{n,1} + f_{n,2})(z_n^d \mu) \\ 
&=& d + (f_{n-2,0} - 2 f_{n-2,0} + f_{n-2,1})(\mu)  \\
& = & d - (f_{n-2,0} - f_{n-2,1})(\mu) \\
& = & d - (d_1 + \dots + d_k) \\
&=& d_0 .
\end{eqnarray*}
Therefore, we have shown that multiplication by $z_n^d$ defines a map from $Q_{n-2}(d_1, \dots, d_k) $ to $Q_{n}(d_0, \dots, d_k)$.  Since multiplication by a fixed element is injective and preserves the partial order, this map an embedding of posets.
\end{proof}

\sse{Remark} Composing the above embedding with the rank-flipping involution on both sides, we get another embedding:
\[ \tau_n z_n^d \tau_{n-2} : Q_{n-2}(d_1, \dots, d_k) \to Q_{n}(d_0, \dots, d_k) \]
which sends $\mu_{n-2}(d_1, \dots, d_k)$ to $\mu_n(d_0, \dots, d_k)$.  We can express this embedding as the composition of the homomorphism $\ka_n:A_{n-2}\to A_n$ which sends each $ z_i \mapsto z_{i+2}$, followed by multiplication by $z_0^d$.  Indeed, for any minimal generator $z_{i_0} \dots z_{i_r}$ of $I_{n-2,r}$ we have:
\[  \tau_n z_n \tau_{n-2}(z_{i_0} \dots z_{i_r}) = \tau_n(z_n z_{n-2-i_0} \dots z_{n-2-i_r}) = z_0 z_{i_0 + 2}  \dots z_{i_r + 2} = z_0 \ka_{n} (z_{i_0} \dots z_{i_r}), \]
from which the general formula follows.

\sse{Remark} The images of the above embeddings can be written down in elementary terms as well:
\[ z_0^d \ka_n Q_{n-2}(d_1, \dots, d_k) = \{ \mu \in Q_n(d_0, \dots, d_k) \mid a_0 = \max_{0 \leq j \leq n-1}(a_j + a_{j+1}) \}, \]
\[ z_n^d Q_{n-2}(d_1, \dots, d_k) = \{ \mu \in Q_n(d_0, \dots, d_k) \mid a_n = \max_{0 \leq j \leq n-1}(a_j + a_{j+1}) \}. \]
These equalities follow immediately from the fact that covering $\{0,1\}$ (resp. $\{n-1,n\}$) will cover the maximum possible amount in each such monomial.

\se{The raising and lowering  algorithm}

The edges of the Hasse diagram of $A_n(m)$ are colored by $\{1, \dots, n \}$, where the $i$-th color corresponds to the following move: 
\[ (a_0, \dots, a_n) \to (a_0, \dots a_{i-1} -1, a_{i} + 1, \dots, a_n).\]
Let $C$ be a saturated chain in $A_n(m)$. We obtain a sequence of colors $(c_1, \dots, c_t)$ by reading the covering relations in $C$ from highest to lowest weight.  We say that a saturated chain $C$ is {\em monotonic} if $c_1 \leq \dots \leq c_t$. 

A pair of consecutive entries $(a_i, a_{i+1})$ of $(a_0, \dots, a_n)$ is called a {\em maximal pair} if:
\[ a_i + a_{i+1} = \max_{0 \leq j \leq n-1} (a_j + a_{j+1}).  \]
Let $\mu = z_0^{a_0}\dots z_n^{a_n} \in Q_n(d_0, \dots, d_k)$.  The following ``right-moving" algorithm produces a monotonic saturated chain between $\mu$ and element of $z_n^d Q_{n-2}(d_1, \dots, d_k)$:

\begin{enumerate}

\item Let $(a_i,a_{i+1})$ be a maximal pair of $(a_0, \dots, a_n)$.  If $i = n-1$, apply the $n$-th color $a_{n-1}$ times and end the chain.

\item If $i < n-1$, apply the $(i+1)$-th color $a_i - a_{i+2}$ times, then go back to the first step and choose the maximal pair $(a_{i+1},a_{i+2})$.

\end{enumerate}

The corresponding ``left-moving" algorithm produces a monotonic saturated chain between $\mu$ and an element of $z_0^d \ka_n Q_{n-2}(d_1, \dots, d_k)$:

\begin{enumerate}

\item Let $(a_{i-1},a_{i})$ be a coverable pair of $(a_0, \dots, a_n)$.  If $i = 1$, apply the inverse of the first color $a_{1}$ times and end the chain.

\item If $i > 1$, apply the inverse of the $i$-th color  $a_i - a_{i-2}$ times, then go back to the first step and choose the maximal pair  $(a_{i-2},a_{i-1})$.

\end{enumerate}

The {\em left transversal chain of} $\mu$ is constructed as follows.  We begin with the leftmost maximal pair in $(a_0,\dots, a_n)$ and run both the left-moving and right moving algorithms.  The concatenation of the two resulting chains is a single monotonic saturated chain passing through $\mu$ which runs between $z_0^d \ka_n Q_{n-2}(d_1, \dots, d_k)$ and $z_n^d Q_{n-2}(d_1, \dots, d_k)$. Similarly, the {\em right transversal chain of} $\mu$ is constructed by applying the same procedure starting with the rightmost maximal pair in $(a_0, \dots, a_n)$.  

\sse{Example} Consider the monomial $\mu = z_0 z_1 z_2 z_4 z_5 \in A_5(5)$.  Since $f_{5,1}(\mu) = 3$ and $f_{5,2}(\mu) = 1$, we see that $\mu \in Q_5(0,1,1)$.  The corresponding lattice vector is $(a_0, \dots, a_5) = (1,1,1,0,1,1)$, which we visualize as 5 indistinguishable balls distributed among 6 labeled urns.  The following diagrams represent the calculation of the left and right transversal chains of $\mu$:
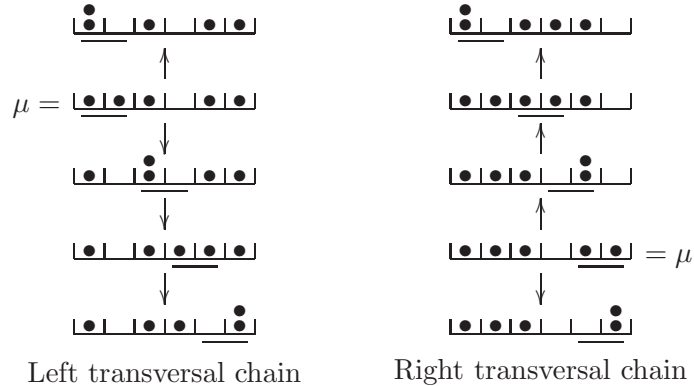
\[ \xy
(0,0); (24,0) **\dir{-}; 
(0,0); (0,2) **\dir{-}; 
(4,0); (4,2) **\dir{-}; 
(8,0); (8,2) **\dir{-}; 
(12,0); (12,2) **\dir{-}; 
(16,0); (16,2) **\dir{-}; 
(20,0); (20,2) **\dir{-}; 
(24,0); (24,2) **\dir{-}; 
(2,1) *+{\bullet}; 
(2,3) *+{\bullet}; 
(10,1) *+{\bullet}; 
(18,1) *+{\bullet}; 
(22,1) *+{\bullet}; 
(1,-1); (7,-1) **\dir{-};
(12,-6); (12,-2) **\dir{-} ?>*\dir{>};
(-5,-10) *+{\mu = };
(0,-10); (24,-10) **\dir{-}; 
(0,-10); (0,-8) **\dir{-}; 
(4,-10); (4,-8) **\dir{-}; 
(8,-10); (8,-8) **\dir{-}; 
(12,-10); (12,-8) **\dir{-}; 
(16,-10); (16,-8) **\dir{-}; 
(20,-10); (20,-8) **\dir{-}; 
(24,-10); (24,-8) **\dir{-}; 
(2,-9) *+{\bullet}; 
(6,-9) *+{\bullet}; 
(10,-9) *+{\bullet}; 
(18,-9) *+{\bullet}; 
(22,-9) *+{\bullet}; 
(1,-11); (7,-11) **\dir{-};
(12,-12); (12,-16) **\dir{-} ?>*\dir{>};
(0,-20); (24,-20) **\dir{-}; 
(0,-20); (0,-18) **\dir{-}; 
(4,-20); (4,-18) **\dir{-}; 
(8,-20); (8,-18) **\dir{-}; 
(12,-20); (12,-18) **\dir{-}; 
(16,-20); (16,-18) **\dir{-}; 
(20,-20); (20,-18) **\dir{-}; 
(24,-20); (24,-18) **\dir{-}; 
(2,-19) *+{\bullet}; 
(10,-19) *+{\bullet}; 
(10,-17) *+{\bullet}; 
(18,-19) *+{\bullet}; 
(22,-19) *+{\bullet}; 
(9,-21); (15,-21) **\dir{-};
(12,-22); (12,-26) **\dir{-} ?>*\dir{>};
(0,-30); (24,-30) **\dir{-}; 
(0,-30); (0,-28) **\dir{-}; 
(4,-30); (4,-28) **\dir{-}; 
(8,-30); (8,-28) **\dir{-}; 
(12,-30); (12,-28) **\dir{-}; 
(16,-30); (16,-28) **\dir{-}; 
(20,-30); (20,-28) **\dir{-}; 
(24,-30); (24,-28) **\dir{-}; 
(2,-29) *+{\bullet}; 
(10,-29) *+{\bullet}; 
(14,-29) *+{\bullet}; 
(18,-29) *+{\bullet}; 
(22,-29) *+{\bullet}; 
(13,-31); (19,-31) **\dir{-};
(12,-32); (12,-36) **\dir{-} ?>*\dir{>};
(0,-40); (24,-40) **\dir{-}; 
(0,-40); (0,-38) **\dir{-}; 
(4,-40); (4,-38) **\dir{-}; 
(8,-40); (8,-38) **\dir{-}; 
(12,-40); (12,-38) **\dir{-}; 
(16,-40); (16,-38) **\dir{-}; 
(20,-40); (20,-38) **\dir{-}; 
(24,-40); (24,-38) **\dir{-}; 
(2,-39) *+{\bullet}; 
(10,-39) *+{\bullet}; 
(14,-39) *+{\bullet}; 
(22,-37) *+{\bullet}; 
(22,-39) *+{\bullet}; 
(17,-41); (23,-41) **\dir{-};
(12,-45) *+{\tn{Left transversal chain}};
(0,0);
(0,0);
(0,0);
(50,0); (74,0) **\dir{-}; 
(50,0); (50,2) **\dir{-}; 
(54,0); (54,2) **\dir{-}; 
(58,0); (58,2) **\dir{-}; 
(62,0); (62,2) **\dir{-}; 
(66,0); (66,2) **\dir{-}; 
(70,0); (70,2) **\dir{-}; 
(74,0); (74,2) **\dir{-}; 
(52,1) *+{\bullet}; 
(52,3) *+{\bullet}; 
(60,1) *+{\bullet}; 
(64,1) *+{\bullet}; 
(68,1) *+{\bullet}; 
(51,-1); (57,-1) **\dir{-};
(62,-6); (62,-2)  **\dir{-} ?>*\dir{>};
(50,-10); (74,-10) **\dir{-}; 
(50,-10); (50,-8) **\dir{-}; 
(54,-10); (54,-8) **\dir{-}; 
(58,-10); (58,-8) **\dir{-}; 
(62,-10); (62,-8) **\dir{-}; 
(66,-10); (66,-8) **\dir{-}; 
(70,-10); (70,-8) **\dir{-}; 
(74,-10); (74,-8) **\dir{-}; 
(52,-9) *+{\bullet}; 
(56,-9) *+{\bullet}; 
(60,-9) *+{\bullet}; 
(64,-9) *+{\bullet}; 
(68,-9) *+{\bullet}; 
(59,-11); (65,-11) **\dir{-};
(62,-16); (62,-12) **\dir{-} ?>*\dir{>};
(50,-20); (74,-20) **\dir{-}; 
(50,-20); (50,-18) **\dir{-}; 
(54,-20); (54,-18) **\dir{-}; 
(58,-20); (58,-18) **\dir{-}; 
(62,-20); (62,-18) **\dir{-}; 
(66,-20); (66,-18) **\dir{-}; 
(70,-20); (70,-18) **\dir{-}; 
(74,-20); (74,-18) **\dir{-}; 
(52,-19) *+{\bullet}; 
(56,-19) *+{\bullet}; 
(60,-19) *+{\bullet}; 
(68,-19) *+{\bullet}; 
(68,-17) *+{\bullet}; 
(63,-21); (69,-21) **\dir{-};
(62,-26); (62,-22) **\dir{-} ?>*\dir{>};
(50,-30); (74,-30) **\dir{-}; 
(50,-30); (50,-28) **\dir{-}; 
(54,-30); (54,-28) **\dir{-}; 
(58,-30); (58,-28) **\dir{-}; 
(62,-30); (62,-28) **\dir{-}; 
(66,-30); (66,-28) **\dir{-}; 
(70,-30); (70,-28) **\dir{-}; 
(74,-30); (74,-28) **\dir{-}; 
(52,-29) *+{\bullet}; 
(56,-29) *+{\bullet}; 
(60,-29) *+{\bullet}; 
(68,-29) *+{\bullet}; 
(72,-29) *+{\bullet}; 
(67,-31); (73,-31) **\dir{-};
(79,-30) *+{= \mu};
(62,-32); (62,-36) **\dir{-} ?>*\dir{>};
(50,-40); (74,-40) **\dir{-}; 
(50,-40); (50,-38) **\dir{-}; 
(54,-40); (54,-38) **\dir{-}; 
(58,-40); (58,-38) **\dir{-}; 
(62,-40); (62,-38) **\dir{-}; 
(66,-40); (66,-38) **\dir{-}; 
(70,-40); (70,-38) **\dir{-}; 
(74,-40); (74,-38) **\dir{-}; 
(52,-39) *+{\bullet}; 
(56,-39) *+{\bullet}; 
(60,-39) *+{\bullet}; 
(72,-37) *+{\bullet}; 
(72,-39) *+{\bullet}; 
(67,-41); (73,-41) **\dir{-};
(62,-45) *+{\tn{Right transversal chain}};
\endxy\]
Upward arrows denote left moves, downward arrows denote right moves, and covered maximal pairs are underlined.  We have omitted those steps where $a_i = a_{i+2}$ because in that case no balls are moved; the only change is the movement of the covering edge.  


\sse{Theorem} Let $n, d_0, \dots, d_k \geq 0$. The transversal chains provide a vertex covering of $Q_n(d_0, \dots, d_k)$ by monotonic saturated chains.  If $d_0 > 0$, then this covering is actually a decomposition.

\begin{proof} For the first statement, the only thing to check is that the algorithm stays within $Q_n(d_0, \dots, d_k)$ at each step.  We claim that if $(a_i , a_{i+1})$ is a maximal pair of $(a_0, \dots, a_n)$ and $a_i > a_{i+2}$, then:
\[  f_{n,r}(a_0, \dots, a_n) = f_{n,r}(a_0, \dots, a_i - 1, a_{i+1} + 1, \dots, a_n) \]
for all $0 \leq r \leq k$.  To prove this, let us deal with each possible case in turn.  Let $\mu'$ denote the monomial obtained by applying the $(i+1)$-th color to $\mu$.  Recall that:
\[ f_{n,r}(\mu) = \min_{F \in \Ga_{n,r}} \al(\ov{F},\mu). \]
If $F \in \Ga_{n,r}$ is a facet such that $\{i,i+1\} \su F$, then $\al(\ov{F},\mu') = \al(\ov{F},\mu)$.  If $i \in F$ and $i+1 \notin F$, then $\al(\ov{F},\mu') = \al(\ov{F},\mu) +1$, which means it is irrelevant for calculating the minimum over such sums.  The only other possibility is that $i \notin F$ and $i+1 \in F$, which means that $i+2 \in F$ as well.  In this case, since $a_i > a_{i+2}$, we see that $\al(\ov{F},\mu) > f_{n,r}(\mu)$, because moving the edge to the left by one step and covering $\{i,i+1\}$ would strictly decrease the sum.  Therefore, 
$ \al(\ov{F},\mu') = \al(\ov{F},\mu) - 1 \geq f_{n,r}(\mu)$, and the overall minimization will be unchanged.


For the second statement, if $d_0 > 0$, then the lexicographically minimal tableau corresponding to $(a_0, \dots, a_n)$ has at least one row of size one.  If the entry in the bottom row is 0 (resp. $n$), then $(a_0, a_1)$ (resp. $(a_{n-1},a_n)$) is the unique maximal pair.  Otherwise, if the entry in the bottom row is $0 < i < n$, then $(a_{i-1},a_i)$ is a maximal pair.  The only other possible maximal pair is $(a_i,a_{i+1})$, which can only happen if $a_{i-1} = a_{i+1}$.  In this case, the algorithm will only move the covering edge, so there is essentially a unique starting maximal pair.  It follows that the left and right transversal chains of $\mu$ coincide and the transversal chains are necessarily disjoint from each other by the uniqueness at each step of the algorithm.  Therefore, if $d_0 > 0$ we obtain a monotonic saturated chain decomposition of $Q_n(d_0, \dots, d_k)$.
\end{proof}

\sse{Corollary}   If $Q_{n-2}(d_1, \dots, d_k)$ has a symmetric chain decomposition and
\[ 1+ 2 d_0 \geq   \sum_{j = 1}^{k} d_j (n-2j)j, \]
then $Q_n(d_0,\dots,d_k)$ has a symmetric chain decomposition.

\begin{proof}
Suppose we have a decomposition:
\[ Q_{n-2}(d_1, \dots, d_k) = C_1 \sqcup \dots \sqcup C_t, \]
where each $C_i$ is a rank-symmetric, saturated chain. We may assume that $C_1$ is the chain containing the unique element $\mu_{n-2}(d_1, \dots, d_k)$ (resp. $\nu_{n-2}(d_1, \dots, d_k)$) of highest (resp. lowest) weight in $Q_{n-2}(d_1, \dots, d_k)$.   Note that the number of edges in any chain of $Q_{n}(d_0, \dots, d_k)$ of maximal length is equal to the weight of $\mu_{n}(d_0, \dots, d_k)$, and:
\begin{eqnarray*}
\tn{wt}(\mu_n(d_0, \dots, d_k)) &=& \sum_{j = 0}^k d_j \tn{wt}(\mu_{n,j}) \\
&=& \sum_{j = 0}^k d_j \sum_{i = 0}^j (n-4i)  \\
&=& \sum_{j = 0}^k d_j (n(j+1) - 2j(j+1)) \\
&=& \sum_{j = 0}^k d_j (n-2j)(j+1)
\end{eqnarray*}
Therefore the number of edges in $C_1$ is equal to:
\[ \sum_{j = 0}^{k-1} d_{j+1} (n-2-2j)(j+1) = \sum_{j = 1}^{k} d_j (n-2j)j. \]
Let $\mu_0$ be a monomial in $C_i$ for some $i$.   
Consider the monomial 
\[ \mu = z_0^{a_0} \dots z_n^{a_n} = \tau_n z_n^d \tau_{n-2} \mu_0 = z_0^d \ka_n \mu_0 \in Q_n(d_0, \dots, d_k).\]
Let us calculate the left transversal chain of $\mu$. The leftmost maximal pair is $(a_0,a_1)$ and $a_1 = 0$ so the first few steps of the algorithm look like:
\[ (a_0, 0 , a_2, a_3, \dots, a_n ) \]
\[ (a_2, a_0 - a_2, a_2, a_3, \dots, a_n) \]
\[ (a_2, a_3, a_0 - a_3 , a_3, \dots, a_n) \]
\[  \vdots \]
\[ (a_2, a_3, \dots, a_n, a_0-a_n ,  a_n) \]
\[ (a_2, a_3, \dots, a_n, 0,  a_0) \]
We see that the $a_0$ term travels to the right while all the other entries shift to the left by two spots. In other words, the left transversal chain of $z_0^d \ka_n\mu_0$ in $Q_{n}(d_0, \dots, d_k)$ has lowest weight monomial $z_n^d \mu_0$. 
The key point is that, for $d_0 > 0$, we obtain a decomposition:
\[ Q_{n}(d_0, \dots, d_k) = R_1 \sqcup \dots \sqcup R_t \]
where $R_i$ is the ``rectangular" poset which contains the chains $z_0^d \ka_n C_i$ and $z_n^d C_i$, as well as all the transversal chains between them.  Furthermore, note that $R_i$ contains $d_0+1$ translates of $z_0^d \ka_n C_i$  (resp. $z_n^d C_i$), namely:
\[ z_0^{d-j} z_1^j \ka_n C_i \;\; (\tn{resp.} \;\; z_n^{d-j} z_{n-1}^j C_i) \;\; \tn{for} \;\; 0 \leq j \leq d_0. \]
The general picture of $R_i$ looks like this:
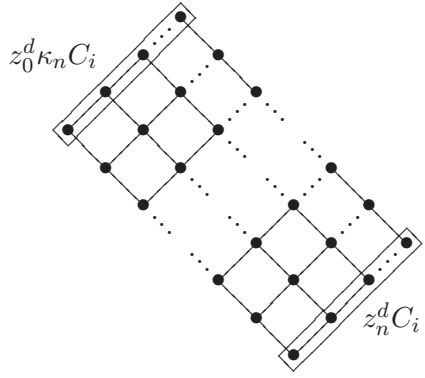
\[\xy
(0,0); (5,5) **\dir{-};
(0,0); (-5,5) **\dir{-};
(5,5); (0,10) **\dir{-};
(5,5); (10,10) **\dir{-};
(-5,5); (0,10) **\dir{-};
(-5,5); (-10,10) **\dir{-};
(11.5,11.5) *+{.};
(12.5,12.5) *+{.};
(13.5,13.5) *+{.};
(0,10); (5,15) **\dir{-};
(10,10); (5,15) **\dir{-};
(15,15); (10,20) **\dir{-};
(6.5,16.5) *+{.};
(7.5,17.5) *+{.};
(8.5,18.5) *+{.};
(0,10); (-5,15) **\dir{-};
(5,15); (0,20) **\dir{-};
(10,20); (5,25) **\dir{-};
(-10,10); (-5,15) **\dir{-}; (0,20) **\dir{-}; 
(1.5,21.5) *+{.};
(2.5,22.5) *+{.};
(3.5,23.5) *+{.};
(0,0) *+{\bullet};
(-5,5) *+{\bullet};
(5,5) *+{\bullet};
(-10,10) *+{\bullet};
(0,10) *+{\bullet};
(10,10) *+{\bullet};
(5,15) *+{\bullet};
(15,15) *+{\bullet};
(10,20) *+{\bullet};
(-5,15) *+{\bullet};
(0,20) *+{\bullet};
(5,25) *+{\bullet};
(-11.5,11.5) *+{.};
(-12.5,12.5) *+{.};
(-13.5,13.5) *+{.};
(-6.5,16.5) *+{.};
(-7.5,17.5) *+{.};
(-8.5,18.5) *+{.};
(-1.5,21.5) *+{.};
(-2.5,22.5) *+{.};
(-3.5,23.5) *+{.};
(3.5,26.5) *+{.};
(2.5,27.5) *+{.};
(1.5,28.5) *+{.};
(-16.5,16.5) *+{.};
(-17.5,17.5) *+{.};
(-18.5,18.5) *+{.};
(-11.5,21.5) *+{.};
(-12.5,22.5) *+{.};
(-13.5,23.5) *+{.};
(-6.5,26.5) *+{.};
(-7.5,27.5) *+{.};
(-8.5,28.5) *+{.};
(-1.5,31.5) *+{.};
(-2.5,32.5) *+{.};
(-3.5,33.5) *+{.};
(-20,20); (-15,25) **\dir{-};
(-20,20); (-25,25) **\dir{-};
(-15,25); (-20,30) **\dir{-};
(-15,25); (-10,30) **\dir{-};
(-25,25); (-20,30) **\dir{-};
(-25,25); (-30,30) **\dir{-};
(-8.5,31.5) *+{.};
(-7.5,32.5) *+{.};
(-6.5,33.5) *+{.};
(-20,30); (-15,35) **\dir{-};
(-10,30); (-15,35) **\dir{-};
(-5,35); (-10,40) **\dir{-};
(-13.5,36.5) *+{.};
(-12.5,37.5) *+{.};
(-11.5,38.5) *+{.};
(-20,30); (-25,35) **\dir{-};
(-15,35); (-20,40) **\dir{-};
(-10,40); (-15,45) **\dir{-};
(-30,30); (-25,35) **\dir{-}; (-20,40) **\dir{-}; 
(-18.5,41.5) *+{.};
(-17.5,42.5) *+{.};
(-16.5,43.5) *+{.};
(-20,20) *+{\bullet};
(-25,25) *+{\bullet};
(-15,25) *+{\bullet};
(-30,30) *+{\bullet};
(-20,30) *+{\bullet};
(-10,30) *+{\bullet};
(-15,35) *+{\bullet};
(-5,35) *+{\bullet};
(-10,40) *+{\bullet};
(-25,35) *+{\bullet};
(-20,40) *+{\bullet};
(-15,45) *+{\bullet};
(0,-2); (17,15) **\dir{-}; (15,17) **\dir{-}; (-2,0) **\dir{-}; (0,-2) **\dir{-}; (13,5) *+{z_n^d C_i};
(-30,28); (-13,45) **\dir{-}; (-15,47) **\dir{-}; (-32,30) **\dir{-}; (-30,28) **\dir{-}; (-32,40) *+{z_0^d \ka_n C_i};
\endxy\]
Now the boundary of this rectangle can be thought of as a disjoint union of a pair of symmetric saturated chains.  If we peel off these two chains, we are left with a smaller rectangle.  We can continue this process as long as we are guaranteed that the boundary of what remains is connected.  At each step, the length of $C_i$ is truncated by two and the length of each transversal chain is truncated by four.  Therefore, for $R_i$ to have a symmetric chain decomposition, it is sufficient that $2d_0+1$ is greater than or equal to the number of edges in a maximal chain of $Q_{n-2}(d_1, \dots, d_k)$. By the above calculation of the length of the chain $C_1$, we conclude that if
\[ 1+ 2 d_0 \geq   \sum_{j = 1}^{k} d_j (n-2j)j, \]
then $Q_n(d_0,\dots,d_k)$ has a symmetric chain decomposition.
\end{proof}

\sse{Remark} The inequality presented in the theorem is not a necessary condition for $Q_n(d_0, \dots, d_k)$ to have a symmetric chain decomposition.  Indeed, as the argument shows, $R_i$ has a symmetric chain decomposition as long as we can successively peel off two boundary chains at a time and leave a rectangle with connected boundary.  The non-transversal edges in these chains do not necessarily have to come from translates of $z_0^d \ka_n C_i$ and $z_n^d C_i$.

\sse{Conjecture} We conjecture that $Q_n(d_0, \dots, d_k)$ is a symmetric chain order for all $n, d_0, \dots, d_{k} \geq 0$.  This statement is strictly stronger than the statement that $L(m,n)$ is a symmetric chain order for all $m,n \geq 0$, since there are symmetric chain decompositions of $L(m,n)$ which do not respect the tropical decomposition.

\sse{Remark}   The set of minimal generators of $I_{n,r}$ is equal to:
\[ \{ z_{i_0} \dots z_{i_r} \mid i_{j+1} - i_j \geq 2 \tn{ for } 0 \leq j \leq r-1 \} = Q_{n}(0, \dots, d_r = 1, \dots, 0). \] 
By mapping $(i_0, \dots, i_r)$ to $(i_0, i_1 -2, \dots, i_r - 2r)$, we see that:
\[ Q_{n}(0, \dots, d_r = 1, \dots, 0) \simeq A_{r+1}(n-2r). \]
Moreover, if we consider monomials whose maximal factorizations involve $d_r$ minimal generators of $I_{n,r}$, we see that:
\[ Q_{n}(0, \dots,d_r, \dots, 0) \simeq A_{r+1}(d_r(n-2r)). \]
In particular, we have that:
\[ Q_2(0,d_1) \simeq A_2(0), \]
which implies that $Q_2(d_0,d_1)$ is a symmetric saturated chain, so $A_2(m)$ is a symmetric chain order for all $m \geq 0$. Now:
\[ Q_3(0,d_1) \simeq A_2(2 d_1), \]
which is a symmetric chain order.  Moreover,
\[ Q_4(0,d_1,0) \simeq A_2(2 d_1)  \quad\tn{and}\quad Q_4(0,0,d_2) \simeq A_3(0), \]
which implies that $Q_4(0,d_0,d_1) \simeq A_2(2d_1)$ is a symmetric chain order.  The cases where $d_0 > 0$ are covered by the theorem, so it follows that $L(m,n)$ is a symmetric chain order for $\min(m,n)\leq 4$. 

\end{document}